\documentclass[12pt]{amsart}
\usepackage{amsmath,amssymb,amsfonts,amsthm}
\usepackage{graphicx}
\usepackage{enumerate}
\usepackage{tikz}
\usepackage[mathscr]{eucal}
\usepackage{float}
\numberwithin{equation}{section}
\newtheorem{theorem}{Theorem}[section]
\newtheorem{definition}{Definition}[section]
\newtheorem{lemma}[theorem]{Lemma}

\newtheorem{note}[theorem]{Note}

\numberwithin{equation}{section}

 \makeatletter
      \def\@setcopyright{}
      \def\serieslogo@{}
      \makeatother

\begin{document}
\bibliographystyle{plain}

\title{Short Paper No. 1: Relation between integer sequence A319433 and the maximum degrees of finite linear Jaco graphs}
\author{Timoth$\acute{e}$e Raso$^{\dag}$ and Johan Kok$^{\ddag}$} 
\address{$^{\dag}$Independent Researcher and Engineer, Lyon, France.}
\email{
timothee.raso@gmail.com}
\address{$^{\ddag}$Independent Mathematics Researcher, City of Tshwane, South Africa.}
\email{jacotype@gmail.com}
\date{}
\begin{abstract}
This short paper presents the interesting relation between integer sequence A319433 in https://oeis.org and the maximum degrees of the family of finite linear Jaco graphs $J_n(x)$ where $n = 1,2,3,\dots$. The result first appeared as a conjecture. Subsequently, a proof through AI-assisted mathematical reasoning came to the fore. Appropriate clarity remarks follow in the main text. It appears that $\textit{Golden ratio}$-like floor function plays an important role in the analysis of the graph structural properties of the family of finite linear Jaco graphs. 
\end{abstract}
\maketitle
\section{Preliminary remarks}
Numerous conjectures appeared in \cite{kok}. Subsequently the first author tested the ability of AI-assisted mathematical reasoning technology (under development) to present valid proofs of the conjectures. The second author evaluated the proofs for correctness. Both authors co-authored this manuscript. See Section \ref{sec1}. 

Most of introduction replicates much which is found in \cite{kok}. The reason is twofold. Firstly, it makes this short paper self-contained and secondly, the arXived manuscript has not officially been published in any journal. This first, in a series of short papers is formally submitted to undergo rigorous peer review as well.
\section{Introduction}
It is assumed that the reader is familiar with the basic terminology and notions of graph theory. See {\cite{bondy, harary, west}. When the context is clear, notation such as $deg_G(v)$ will be abbreviated to $deg(v)$ and so on. In 2014 the authors Kok et al. introduced the family of directed finite linear Jaco graphs. See \cite{kok1, kok2} and the references thereto. Unless stated otherwise, all graphs in this short paper will be the underlying graph of the directed finite linear Jaco graphs (for reference, linear Jaco graphs) . 

The motivation for the study in \cite{kok} does not lie in possible applications of linear Jaco graphs \textit{per se} but rather in a new research avenue i.e. relations between some integer sequences and certain parameters of linear Jaco graphs. In addition to the aforesaid the validation of AI-assisted mathematical reasoning to prove of such conjectured relations has come to the fore. Clearly, the aforesaid has real world applications.

Recall that a sequence of non-negative integers, say
\begin{center}
$\textbf{d} = (d_1,d_2,d_3,\dots,d_n)$
\end{center}
is said to be graphical if there exists a simple graph $G$ which has degree sequence \textbf{d}. Put differently, \textbf{d} is graphical in respect of the vertex degrees if there exists a simple graph $G$ which has degree sequence \textbf{d}. Note that \textbf{d} could be either a random or a definable sequence. A definable sequence can only be generated from a closed or recursive formula. A generalization can be investigated. That is, if $p(G)$ is a graph parameter (independence number, domination number, chromatic number, hub number and so on) then a definable sequence (finite or infinite) $\textbf{p} = p_1,p_2,p_3,\dots,p_n$ is $p$-graphical if a family (or set) or subset of a family of simple graphs exists, say $\mathcal{F} = \{G_1, G_2, G_3,\dots,G_n\}$ such that,
\begin{center}
$p(G_1) = p_1$, $p(G_2) = p_2$, $p(G_3) = p_3,\cdots,p(G_n) = p_n$.
\end{center}
Note that the notion of a $p$-graphical sequence extends over graph parameters (global as mentioned above) and parameters related to graph elements (local), for example the eccentricity of a vertex, $e(v)$, the degree of a vertex, $deg_G(v)$, the distance between two distinct edges and so on. For the aforesaid the notation is $p(G,v)$ or $p(G,e)$. The parameter may also be a hybrid such as vertex-edge domination and edge-vertex domination (see \cite{peters}) or it may concern subsets of $V(G)$ or $E(G)$ and so on. The corresponding notation will be $p(G,(v,e))$ or $p(G,(e,v))$ or $p(G,S)$ where $S \subseteq V(G)$ or $S \subseteq E(G)$.

\textbf{Example:} The infinite integer sequence $\textbf{p} = 1,2,3,\dots$ is $p$-graphical for $ p = \chi(G)$ (chromatic number of graph) since $\mathcal{F} = \{K_1, K_2, K_3,\dots\}$ exists and
\begin{center}
$\chi(K_1) = 1$, $\chi(K_2) = 2$, $\chi(K_3) = 3,\cdots$
\end{center}
\begin{note}
The notion of graphical degree sequences is generalized to \textit{graphical graph parametric sequences}. This avenue for further research remains open. 
\end{note}
\section{On linear Jaco graphs}
Readers are referred to \cite{kok1, kok2} for the formal definitions of $J_\infty(x)$ and $J_n(x)$. Since our interest lies with the underlying graph we recall the construction method as explained in \cite{kok3} with some additional refinement.  

\textit{Construction method to obtain $J_n(x)$, $n \geq 4$:} 

(a) Let $J_1(x) = P_1$, $J_2(x) = P_2$ and $J_3(x) = P_3$.

(b) Let the empty graph (edgeless graph) of order $m \geq 1$ be denoted by $\mathfrak{N}_m$. Let $X = \{v_i: i = 4,5,6,\dots,n\}$

and let $G_2 = P_3 \cup \mathfrak{N}_{n-3}$ where $V(\mathfrak{N}_{n-3}) = X$.

(c) By normal consecutive step-count for $i = 3,4,5,\dots, n-1$ do as follows: 

To obtain $G_i$ add the edges $v_iv_{i+1}$, $v_iv_{i+2}$,$\dots, v_iv_{i+t}$ with $t$ a maximum such that, $deg_{G_i}(v_i) \leq i$.

(d) After completion of step-count $i = n-1$ label the resultant graph $J_n(x)$.
Note that $J_n(x)$ is a simple connected undirected graph.

Recall from \cite{kok3} that for a given $J_n(x)$ the open neighborhood of $v_i$ in a given by:
$N(v_i) = N^-(v_i) \cup N^+(v_i)$ where $N^-(v_i) = \{v_j:j < i$ and the edge $v_jv_i$ exists$\}$ and $N^+(v_i) = \{v_j:j > i$ and the edge $v_iv_j$ exists$\}$.
Also let $t_1(v_i) = |N^-(v_i)|$ and $t_2(v_i) = |N^+(v_i)|$ and let $t^\star_2(v_i) = max\{t_2(v_i)$ in $J_n(x)$ evaluated at $v_i$: $n \geq i\}$. The aforesaid yields $t^\star_2(v_i) = i - t_1(v_i)$.

All sequences mentioned in this paper are found in \cite{sloane}. Recall that Wagner \cite{kok1, kok3} identified that in accordance with sequence A060144 the result $t_1(v_i) = \lfloor \frac{2(i+1)}{3+\sqrt{5}}\rfloor$, $i = 1,2,3,\dots$ holds. Put differently, the sequence A060144 is $p$-graphical where $p(G) = t_1(v_i)$, $i = 1,2,3,\dots$ and $\mathcal{F} = \{G:G=J_n(x), n = 1,2,3,\dots\}$ or simply $J_\infty(x)$.

\subsection{Maximum degree of $J_n(x)$, $n = 1,2,3,\dots$:} 
In \cite{kok} it is conjectured (see Conjecture 2.2 in \cite{kok}) that for the finite linear Jaco graphs $J_n(x)$, $ n = 2,3,4,\dots$ the respective maximum degrees $\Delta(J_n(x))$ are yielded by the sequence A319433 with closed formula $\Delta(J_n(x)) = \lfloor \frac{2(n+2)}{1+ \sqrt{5}}\rfloor -1$, $n = 2,3,4,\dots$

Furthermore, it implies that sequence A319433 is $p$-graphical where $p(G) = \Delta(G)$ and $\mathcal{F} = \{G:G=J_n(x), n = 2,3,\dots\}$.

The AI-assisted mathematical reasoning included $n=1$ and it is easy to see why. Since $1 < \frac{2(1+2)}{1+ \sqrt{5}} <2$ it follows that $\lfloor \frac{2(1+2)}{1+ \sqrt{5}}\rfloor -1 = 0$. Therefore, the conjecture may be considered for $n = 1,2,3,\dots$ The AI-assisted mathematical reasoning utilized the conventional directed linear Jaco graph with the notation thereof as found in the founding papers. See \cite{kok1, kok2}. Adapting the proof to the notation in \cite{kok} is easy in that:
\begin{center}
$t_1(v_i) = d^-(v_i)$, $t^*_2(v_i) = d^+(v_i)$ in $J_\infty(x)$ and

$t_2(v_i) = d^+(v_i)$ in the finite truncation to obtain $J_n(x)$, $n \in \mathbb{N}$. 
\end{center}
Note that for a sufficiently large $n$ (which always exists), $t^*_2(v_i) = t_2(v_i) = d^+(v_i)$.
\section{AI-assisted mathematical reasoning:} It is known that the golden ratio is,
\begin{center}
$\varphi = \frac{1+\sqrt{5}}{2}$.
\end{center}
The known identities which will be used frequently are,
\begin{center}
$\varphi^2 = \varphi +1$, $\frac{1}{\varphi} = \varphi-1$, $\frac{1}{\varphi^2} = 2-\varphi$, $\frac{1}{\varphi} + \frac{1}{\varphi^2} = 1.$
\end{center}
Note that since $\varphi$ is irrational, every nonzero integer multiple of $\varphi$, $\frac{1}{\varphi}$, $\frac{1}{\varphi^2}$ is irrational. Also, for a given $v_i$ the largest $j>i$ such that $v_j$ is a neighbor of $v_i$ is called the \emph{reach} of $v_i$ and denoted by $R_i$ and defined by $R_i = 2i - t_1(v_i)$. It implies that a vertex $v_j$, $i < j \leq R_i$ is a neighbor of $v_i$. Consequently,
\begin{center}
$t^*_2(v_i) = R_i -i = i-t_1(v_i)$.
\end{center}
In a finite linear Jaco graph the degree of a vertex $v_i$ is bounded by $deg(v_i) \leq i$. For use later on the next is an important definition.
\begin{definition}
For $i \geq 1$, define
\begin{center}
$A_i = \lfloor \frac{i+1}{\varphi^2}\rfloor$ and $B_i = \lfloor \frac{i+1}{\varphi}\rfloor$.
\end{center}
\end{definition}
\begin{lemma}\label{lem1}
For integers $i \geq 1$ we have $A_i + B_i = i$.
\end{lemma}
\begin{proof}
Note that both $\frac{i+1}{\varphi^2}$ and $\frac{i+1}{\varphi}$ are non-integers. However, from the fourth identity i.e. $\frac{1}{\varphi} + \frac{1}{\varphi^2} = 1$ it follows that, 
\begin{center}
$\frac{i+1}{\varphi^2} + \frac{i+1}{\varphi} = i+1$ and an integer.
\end{center}
Hence, the identity $\lfloor x\rfloor + \lfloor y\rfloor = x+y-1$ applies. Therefore,
\begin{center}
$A_i + B_i = \lfloor \frac{i+1}{\varphi^2}\rfloor + \lfloor \frac{i+1}{\varphi}\rfloor = \underbrace{\frac{i+1}{\varphi^2} + \frac{i+1}{\varphi}}_{= i+1} -1 = i.$
\end{center} 
\end{proof}
\begin{lemma}(Reach identity)\label{lem2}
For every $i \geq 1$ we have, $i+B_i = \lfloor (i+1)\varphi\rfloor -1$ and in particular, the sequence $i+B_i$ is strictly increasing.
\end{lemma}
\begin{proof}
From the second identity i.e. $\varphi = 1 + \frac{i}{\varphi}$ it follows that,
\begin{center}
$(i+1)\varphi = i+1+ \frac {i+1}{\varphi}$. 
\end{center}
Since $i+1$ is an integer it follows that, 
\begin{center}
$\lfloor (i+1)\varphi\rfloor = i+1 + \lfloor \frac {i+1}{\varphi}\rfloor = i+1+B_i$.
\end{center}
This settles the result.

Also, the claim to strict increase follows because $1 <\varphi < 2$ hence, $\lfloor (i+1)\varphi\rfloor$ increases by $1~or~2$.
\end{proof}
\subsection{Deriving $t_1(v_i)$ and $t^*_2(v_i)$}
The AI-assisted mathematical reasoning included autonomous proofs of two important results.
\begin{theorem}\label{thm1}
For every $i \geq 1$ we have, $t_1(v_i) = A_i = \lfloor \frac{i+1}{\varphi^2}\rfloor.$ Consequently, in $J_\infty(x)$ we have that, $t^*_2(v_i) = B_i = \lfloor \frac{i+1}{\varphi}\rfloor.$
\end{theorem}
\begin{proof}
The proof will be by mathematical induction.

For $i=1$, there is no smaller indexed vertex, hence $t_1(v_1) = 0$. Also $A_1 = \lfloor \frac{1+1}{\varphi^2}\rfloor = 0$. Thus the result holds for $i=1$. Assume the results holds for all $i = 1,2,3,\dots, j-1$, $j \geq 2$. For each $i < j$ the induction hypothesis read with Lemma \ref{lem1} give
\begin{center}
$R_i = 2i -t_1(v_i) = 2i - A_i = i + (i-A_i) = i+B_i$.
\end{center}
By Lemma \ref{lem2}, $R_i = \lfloor (i+1)\varphi\rfloor -1$. The adjacency rule permits the edge $v_iv_j$, $i<j$ when $R_i \geq j$. Since the sequence $R_i$ is strictly increasing, there exists a terminal interval $m, m+1, m+2,\dots,j-1$ or put differently, a subset of vertices $\{v_m, v_{m+1}, v_{m+2},\dots, v_{j-1}\}$ which permits the edges $v_iv_j$, $i = m, m+1, m+2, \dots, j-1$. Furthermore, $m$ is the smallest integer satisfying $R_m \geq j$. Also, $R_m \geq j \Longleftrightarrow \lfloor (m+1)\varphi\rfloor \geq j+1$. Since $(m+1)\varphi$ is irrational it follows that, 
\begin{center}
$\lfloor (m+1)\varphi\rfloor \geq j+1 \Longleftrightarrow (m+1)\varphi > j+1$.
\end{center}
Since $\frac{j+1}{\varphi}$ is non-integral, the smallest integer strictly greater than $\frac{j+1}{\varphi}$ is, $m+1 > \frac{j+1}{\varphi}$.  Therefore, the largest integer strictly less than $\frac{j+1}{\varphi}$ is, $m = \lfloor \frac{j+1}{\varphi}\rfloor = B_j$. So it follows that the number of vertices $v_i$, $i <j$ such that the edges $v_iv_j$ exist is given by,
\begin{center}
$t_1(v_j) = j-m = j- B_j$.
\end{center}
By Lemma \ref{lem1} $j - B_j = A_j$ or put differently, $t_1(v_j) = A_j$. This settles the induction proof of the first part of the theorem.

Finally the fact that $t^*_2(v_i) = i-t_1(v_i)$ read together with Lemma \ref{lem1} concludes that $t^*_2(v_i) = i - A_i = B_i$ which proves the second part of the theorem.
\end{proof}
\begin{note}
Theorem \ref{thm1} independently confirms what Wagner \cite{kok1, kok3} established i.e. $t_1(v_i) = \lfloor \frac{2(i+1)}{3+\sqrt{5}}\rfloor$ because $3+\sqrt{5} = 2\varphi^2$.
\end{note}
\subsection{Deriving degree of vertex $v_i$ in $J_n(x)$, $n \in \mathbb{N}$, $1 \leq i \leq n$} 
Note that for a given $i \in \mathbb{N}$, $t_1(v_i)$ remains a constant in the family of finite linear Jaco graphs $J_n(x)$, $n = 1,2,3,\dots$. However, due to truncation to obtain a finite linear Jaco graph $J_n(x)$, $n \in \mathbb{N}$ it follows that $t_2(v_i) \leq t^*_2(v_i)$. The next theorem provides a result for $deg(v_i)$, $v_i \in V(J_n(x))$, $n \in \mathbb{N}$.
\begin{theorem}\label{thm2}
For every $n \geq 1$ and every $1 \leq i \leq n$ and every $v_i \in V(J_n(x))$,
\begin{center}
$deg(v_i) = min\{i, n-\lfloor \frac{i+1}{\varphi}\rfloor\}.$
\end{center}
\end{theorem}
\begin{proof}
By Theorem \ref{thm1} the vertex $v_i$ has $t_1(v_i) = A_i$ and in $J_\infty(x)$ we have $t^*_2(v_i) = B_i$. The neighbors corresponding to $t^*_2(v_i)$ are,
\begin{center}
$v_{i+1}, v_{i+2}, v_{i+3}, \dots, v_{i+B_i}$.
\end{center}
In $J_n(x)$ the neighbors corresponding to $t_2(v_i)$ is limited to a subset of the set $\{v_{i+1}, v_{i+2}, v_{i+3}, \dots, v_n\}$. Hence the number (cardinality of permissible subset) of such neighbors is $min\{B_i, n-i\}$. Hence, $deg(v_i) = A_i + min\{B_i, n-i\}$. Using $A_i + B_I =i$ it follows that,
\begin{center}
$deg(v_i) = min\{A_i+ B_i, A_i + n-i\}$

$= min\{i, n-(i-A_i)\}$

$= min\{i, n-B_i\}$.
\end{center}
By substituting $B_i = \lfloor \frac{i+1}{\varphi}\rfloor$ the result is settled.
\end{proof}
\begin{note}
No boundary approximation is involved in Theorem \ref{thm2}. Truncation removes precisely those neighbors whose indices exceed $n$. The indices of the forward neighbors corresponding to $t_2(v_i)$ form a consecutive interval of integers and the number of such integers is exactly $min\{B_i, n-i\}$.
\end{note}
\subsection{A scalar threshold for $\Delta(J_n(x))$} Define $\Delta_n = \Delta(J_n(x))$. By Theorem \ref{thm2} we have $\Delta_n = \underset{1\leq i\leq n}{max}\{min\{i, n- B_i\}\}$. For $i = 1,2,3,\dots$ the sequence $B_i$ is non-decreasing, since it is the floor of a strictly increasing function.
\begin{lemma}(Threshold characterization)\label{lem3}
For every $n \geq 2$ and every integer $r$, $1 \leq r \leq n$,
\begin{center}
$\Delta_n \geq r \Leftrightarrow r+B_r \leq n$.
\end{center}
\end{lemma}
\begin{proof}
From $\Delta_n = \underset{1\leq i\leq n}{max}\{min\{i, n- B_i\}\}$ it follows that $\Delta_n \geq r$ if and only if there exists an $i$, $1 \leq i \leq n$ such that both $i \geq r$ and $n-B_i \geq r$. Equivalently, there exists an $i \in \{r, r+1, r+2,\dots, n\}$ such that $B_i \leq n-r$. 

Since $B_i$ is non-decreasing there exists an index $i \in \{r, r+1, r+2,\dots, n\}$ satisfying
\begin{center}
$B_i \leq n-r$ if and only if $B_r \leq n-r$. 
\end{center}
Therefore, $\Delta_n \geq r$ if and only if $r + B_r \leq  n$.
\end{proof}
Lemma \ref{lem3} is the decisive simplification in that, graph theoretical maximization is replaced by inversion of the one-variable quantity $r + B_r$.
\subsection{Proof of main conjecture}
\begin{theorem}(Maximum degree of each $J_n(x)$, $n = 1,2,3,\dots$)\label{thm3}
For every integer $n \geq 1$ we have,
\begin{center}
$\Delta(J_n(x)) = \lfloor \frac{n+2}{\varphi}\rfloor -1 = \lfloor \frac{2(n+2)}{1 + \sqrt{5}}\rfloor -1.$
\end{center}
\end{theorem}
\begin{proof}
Case 1: For $n=1$ the result is trivial because $\Delta(J_1(x))=0$.

Case 2: For $n \geq 2$: By Lemma \ref{lem3} for every $r$ with $1 \leq r \leq n$ it is true that,
\begin{center}
$\Delta_n \geq r \Leftrightarrow r + \lfloor \frac{r+1}{\varphi}\rfloor \leq n$

$\Leftrightarrow \lfloor (r+1)\varphi\rfloor \leq n+1$

$\Leftrightarrow (r+1)\varphi < n+2$

$\Leftrightarrow r+1 < \frac{n+2}{\varphi}$.
\end{center}
The number $\frac{n+2}{\varphi}$ is irrational hence, the largest integer $r+1$ below $\frac{n+2}{\varphi}$ is $\lfloor \frac{n+2}{\varphi}\rfloor$. Consequently, the largest integer $r$ for which $\Delta_n \geq r$ is  $r = \lfloor \frac{n+2}{\varphi}\rfloor -1$. Since $\frac{1}{\varphi} = \frac{2}{1+\sqrt{5}}$ the result is settled.
\end{proof}
\section{Conclusion}
Theorem \ref{thm3} essentially proves the Conjecture 2.2 in \cite{kok}. It can be said that integer sequence A319433 found in \cite{sloane} is $p$-graphical where $p(G) = \Delta(G)$ and $\mathcal{F} = \{G:G=J_n(x), n = 1,2,3,\dots\}$.

\section{AI disclosure}\label{sec1} 
AI-assisted mathematical reasoning technology was used substantially in the mathematical exploration, development of candidate arguments, intermediate verification, and preparation of the source manuscript. The technology was used within an author-directed problem-solving and verification methodology. All resulting mathematical claims and proofs were subsequently reviewed and validated before human drafting of a final manuscript. The human authors take full responsibility for all mathematical statements, proofs, computations, citations, and conclusions presented herein.

After completion of this captioned topic and others found in \cite{kok} the first author became aware of an independent AI-proof of Conjecture 2.2 posted by Shivan Patel on the repository MathDB. Accordingly no claim of priority is made here.
\section*{Dedication}
This short paper is dedicated to late Mr Fransloo Vorster (nicknamed "Fifi") who was the teacher of mathematics of the second author at Bellville High School, Western Cape Province, South Africa. The dedication serves to acknowledge the remarkable way in which Mr Vorster could explain the wonders of High School Mathematics. May his soul rest in peace.


\begin{thebibliography}{999}
\bibitem{bondy} J.A. Bondy and U.S.R. Murty, {\it Graph Theory with Applications}, Macmillan Press, London, (1976).  
\bibitem{harary} F. Harary, {\it Graph Theory}, Addison-Wesley, Reading MA, (1969).
\bibitem{kok} J. Kok, Integer sequences with conjectured relation with certain graph parameters of the family of linear Jaco graphs, {\it arXiv:2507.16500} (2025).
\bibitem{kok1} J. Kok, S. Naduvath, C. Kalthavalappil, G.K. Augusthy and U.M. Augustine, Special issue on: Reflections on Jaco Graphs, {\it Journal of Informatics and Mathematical Sciences}, {\bf 8}(2), (2016).
\bibitem{kok2} J. Kok, C. Susanth and S.J. Kalayathankal, A study on linear Jaco graphs, {\it Journal of Informatics and Mathematical Sciences}, {\bf 7}(2), 69--80,(2015).
\bibitem{kok3} J. Kok, Research note: Domination of exact deg-centric Jaco graphs, {\it Open Journal of Discrete Applied Mathematics}, {\bf 8}(2), 52--60,(2025).
\bibitem{patel} S. Patel, Maximum-degree formula conjecture for finite linear Jaco graphs, {\it MathDB, August} (2026).
\bibitem{peters} K.W. Peters, {\it Theoretical and Algorithmic Results on Domination and Connectivity}, Ph.D. thesis, Clemson University, Clemson SC, (1986).
\bibitem{sloane} N.J.A. Sloane et al., {\it The On-Line Encyclopedia of Integer Sequences}, https://oeis.org.
\bibitem{west} D.B. West, {\it Introduction to graph theory}, Prentice Hall of India, New Delhi, (2001).
\end{thebibliography}
\end{document}